\documentclass[a4paper,11pt]{article}

\usepackage[utf8]{inputenc}
\usepackage[T2A]{fontenc}
\usepackage[russian,english]{babel}

\usepackage{graphicx}
\usepackage{amsfonts}
\usepackage{amsthm}
\usepackage{amsmath}
\usepackage{amssymb}

\usepackage{accents}
\usepackage{indentfirst}
\usepackage{float}
\usepackage{mathtools}
\usepackage[noadjust]{cite}
\usepackage{hyperref}

\usepackage{tikz}

\usepackage{enumitem}
\setlist[enumerate]{font=\upshape,label=(\arabic*)}
\setlist{noitemsep,topsep=4pt} 

\newtheorem{theorem}{Theorem}
\newtheorem{proposition}[theorem]{Proposition}
\newtheorem{corollary}[theorem]{Corollary}
\newtheorem{lemma}[theorem]{Lemma}
\theoremstyle{definition}
\newtheorem{definition}[theorem]{Definition}

\theoremstyle{definition}

\newtheorem{example}[theorem]{Example}

\numberwithin{theorem}{section}
\numberwithin{equation}{section} 
\numberwithin{figure}{section} 

\usepackage{geometry}

\providecommand{\keywords}[1]{\textbf{Keywords:} #1}
\providecommand{\msc}[1]{\textbf{MSC 2020:} #1}

\newcommand\freefootnote[1]{%
    \bgroup
    \renewcommand\thefootnote{\fnsymbol{footnote}}%
    \renewcommand\thempfootnote{\fnsymbol{mpfootnote}}%
    \footnotetext[0]{#1}%
    \egroup
}

\begin{document}
\title{\vspace{-0.2cm}Additive mappings of unweighted directed graphs preserving individual values of the cycle length}

\author{
Alexander Vlasov$^{a,b}$, Alexander Guterman$^{c}$, and Svetlana Zhilina$^{a,b}$
}
\date{\small \em
$^a$Lomonosov Moscow State University, Moscow, 119991, Russia\\
$^b$Moscow Center of Fundamental and Applied Mathematics, Moscow, 119991, Russia\\
$^c$Bar-Ilan University, Ramat Gan, 5290002, Israel
}

\freefootnote{The work of the first and the third authors was supported by the Ministry of Science and Higher Education of the Russian Federation as part of the program of the Moscow Center of Fundamental and Applied Mathematics of Lomonosov Moscow State University under agreement No. 075-15-2025-345.}

\freefootnote{Email addresses:  \texttt{aleksandr.vlasov@math.msu.ru} (Alexander Vlasov), \texttt{alexander.guterman@biu.ac.il} (Alexander Guterman), \texttt{s.a.zhilina@gmail.com} (Svetlana Zhilina)}

\maketitle

\begin{abstract}
On the set of unweighted directed graphs on $n$ vertices, the union of graphs is defined as the union of their arc sets. A mapping of graphs is said to be additive if it preserves the union operation. In this paper, a characterization of additive mappings of unweighted directed graphs preserving a fixed value $k \geq 3$ of the cycle length is obtained.
\end{abstract}


\msc{05C38, 05C50, 15A80, 15A86.}

\section{Introduction}

A set $\mathcal{S}$ with two binary operations, addition and multiplication, is called a \textit{semiring} if
\begin{enumerate}
    \item $\mathcal{S}$ is a commutative monoid under addition with identity element $0$;
    \item $\mathcal{S}$ is a semigroup under multiplication;
    \item multiplication is distributive over addition;
    \item $a\cdot 0 = 0\cdot a = 0$ for any $a\in \mathcal{S}$.
\end{enumerate}

The \textit{tropical semiring} $\mathbb{R}_{\max} = \mathbb{R}\cup \{-\infty\}$ plays a special role in algebra. The operations of tropical addition $\oplus$ and tropical multiplication $\otimes$ are defined on it as follows:
\begin{align*}
    &a\oplus b = \max\{a, b\},\\
    &a\otimes b = a+b
    & &\text{for all }a, b\in \mathbb{R},\\
    &a\oplus (-\infty) = (-\infty)\oplus a = a,\\
    &a\otimes (-\infty) = (-\infty)\otimes a = -\infty
    & &\text{for any } a\in \mathbb{R}_{\max}.
\end{align*}
The element $-\infty$ is the identity element with respect to tropical addition, and $0$ is the identity element with respect to tropical multiplication.

If $\mathcal{S}$ is a semiring, then the set $M_n(\mathcal{S})$ of $n\times n$ matrices over $\mathcal{S}$ is also a semiring. For brevity, we denote the set of matrices over the tropical semiring by~$M_n$. Throughout the paper we assume that $n\geq 2$.

When studying tropical matrices, it is convenient to consider weighted directed graphs. Namely, to each matrix $A\in M_n$ we associate the graph $G(A)$ on the vertices $1, \ldots, n$ such that the arc from vertex $i$ to vertex $j$ has weight $A_{ij}$ if $A_{ij}\neq -\infty$, and is absent otherwise. Moreover, to a matrix $A\in M_n$ we can also associate the unweighted directed graph $\widehat G(A)$ obtained from $G(A)$ by deleting the weights of all arcs. All graphs under consideration have no multiple arcs but may contain loops. A matrix $A$ is called \textit{irreducible} if the graph $G(A)$ is strongly connected. By a \textit{cycle} we always mean a simple directed cycle, that is, a closed path with no repeated vertices except that the initial and the terminal vertices coincide. A loop is a cycle of length~$1$.

The correspondence between tropical matrices and directed graphs finds applications in many theoretical and practical problems. For instance, tropical matrices play an important role in graph algorithms, scheduling theory, and the analysis of manufacturing and logistics systems. They are also used in data analysis and machine learning: for example, tropical matrix factorization is used for the recovery of missing values in sparse data.

One of the key invariants of tropical matrices is the \textit{cyclicity index}. It is defined first for directed graphs and then for tropical matrices.

\begin{definition} \label{definition:cyclicity-index}
    Let $G$ be a directed graph. Its \textit{cyclicity index} $\sigma_G$ is defined as follows:
    \begin{itemize}
        \item if $G$ is strongly connected and contains no cycles, then $\sigma_G = 1$;
        \item if $G$ is strongly connected and contains at least one cycle, then $\sigma_G$ is equal to the greatest common divisor of the lengths of all cycles in~$G$;
        \item if $G$ is not strongly connected, then $\sigma_G$ is equal to the least common multiple of the cyclicity indices of all its strongly connected components.
    \end{itemize}
\end{definition}

The cyclicity index of a tropical matrix $A$ is defined as the cyclicity index of the critical subgraph $G_c(A)$ of the graph~$G(A)$, that is, of the union of all cycles with the largest average weight in~$G(A)$. If $G(A)$ contains no cycles, we set $\sigma(A) = 1$. As the following theorem shows, the cyclicity index plays an important role in the study of the behavior of powers of tropical matrices.

\begin{theorem}[{\cite[Theorem~3.9]{HOW}}]
    Let $A\in M_n$ be an irreducible tropical matrix with finite eigenvalue $\lambda$ and cyclicity index $\sigma$. Then there exists a positive integer~$N$ such that for any positive integer $t\geq N$ we have
    \[
        A^{\otimes(t+\sigma)} = \lambda^{\otimes \sigma}\otimes A^{\otimes t}.
    \]
\end{theorem}

Once an invariant is defined, the question naturally arises as to which transformations preserve it. One of the first to study a problem of this kind was Frobenius~\cite{F}, who characterized the bijective linear mappings of complex matrices that preserve the determinant. Later this result was substantially improved by Schur~\cite{Schur} and Dieudonn\'e~\cite{Di}. Subsequently, many authors investigated mappings that preserve various invariants of matrices over various algebraic structures; see~\cite{BP1,BP2,GM1,Pierce}  and the references therein.  

We say that a mapping $T: M_n\to M_n$ \textit{preserves the cyclicity index} if $\sigma(T(A)) = \sigma(A)$ for any matrix $A\in M_n$. A mapping $T$ \textit{preserves the value $k$ of the cyclicity index} if for any matrix $A\in M_n$ the equality $\sigma(A) = k$ implies $\sigma(T(A)) = k$. A mapping $T$ is called \textit{linear} if
\[
    T(\alpha\otimes A\oplus \beta\otimes B)
    = \alpha\otimes T(A)\oplus \beta\otimes T(B)
\]
for any $A, B\in M_n$ and $\alpha, \beta\in \mathbb{R}_{\max}$.

Linear mappings of tropical matrices that preserve the cyclicity index were previously considered in~\cite{GKT,GKV1,GKV2}. In~\cite{GKT} a description of bijective linear mappings preserving the cyclicity index was obtained, and in~\cite{GKV1} a complete characterization of such mappings without the bijectivity assumption was given. In~\cite{GKV2} it was shown that preserving all values of the cyclicity index is equivalent to preserving its values $1$ and~$2$ simultaneously. In the subsequent paper~\cite{GVZ2} we consider bijective linear mappings of tropical matrices preserving individual values of the cyclicity index. In particular, we obtain a characterization of the mappings preserving the value~$1$, and investigate the case of preserving one fixed value $k\geq 3$ for $n\geq \max(k, 4)$.

Let $\widehat M_n$ be the set of unweighted directed graphs without multiple arcs on the fixed vertex set $\{1, \ldots, n\}$. We define the union of graphs as the union of their arc sets. We call a mapping $S: \widehat M_n\to \widehat M_n$ \textit{additive} if
\[
    S(\widehat G_1\cup \widehat G_2)
    = S(\widehat G_1)\cup S(\widehat G_2)
\]
for any $\widehat G_1, \widehat G_2\in \widehat M_n$. It is clear that for any matrices $A, B \in M_n$ and $\alpha \in \mathbb{R}$ we have $\widehat G(A \oplus B) = \widehat G(A) \cup \widehat G(B)$ and $\widehat G(\alpha \otimes A) = \widehat G(A)$. Hence every linear mapping of matrices $T: M_n\to M_n$ induces an additive mapping of graphs $\widehat{T}: \widehat M_n \to \widehat M_n$ given by the formula $\widehat{T}(\widehat G(A)) = \widehat G(T(A))$ for all $A \in M_n$. Indeed, it suffices to define the values of $\widehat{T}$ on the graphs containing only one arc and then extend $\widehat{T}$ by additivity. Moreover, if $T$ is bijective, then $\widehat{T}$ is also bijective.

\begin{definition}\label{def_cycle_length_preserving}
    We say that a mapping $S: \widehat M_n\to \widehat M_n$ \textit{preserves the value $k$ of the cycle length} if it maps every graph whose arcs are exactly the arcs of a single cycle of length~$k$ to a graph of the same form. A mapping $S$ \textit{preserves cycle lengths} if it preserves each value of the cycle length.
\end{definition}

In~\cite{GVZ2} we show that:
\begin{enumerate}
    \item if a bijective linear mapping $T: M_n\to M_n$ preserves the value $1$ of the cyclicity index, then the induced mapping $\widehat T:~\widehat M_n\to\widehat M_n$ preserves cycle lengths;
    \item if a bijective linear mapping $T: M_n\to M_n$ preserves the value $k$ of the cyclicity index, $k\leq n$, then the induced mapping $\widehat T$ preserves the value $k$ of the cycle length.
\end{enumerate}

Thus, in order to describe linear mappings of tropical matrices preserving individual values of the cyclicity index, it is necessary to characterize additive mappings of unweighted directed graphs preserving individual values of the cycle length. The present paper is devoted to this problem.

Note that this topic is also of independent interest. For instance, Thomassen~\cite{Thomassen} proved that if there exists a bijective mapping between the arcs of two strongly connected tournaments with finite vertex sets such that both this mapping and its inverse preserve the value $4$ of the cycle length, then these tournaments are either isomorphic or anti-isomorphic. Furthermore, Woodall~\cite{Woodall} showed that if $G$ and $H$ are undirected graphs without cut-edges, then a length-preserving isomorphism between the integer cycle groups of $G$ and $H$ exists if and only if $G$ and $H$ have isomorphic circuit matroids.

We remark that the mapping $S: \widehat M_n\to \widehat M_n$ preserves the value $1$ of the cycle length if and only if the set of loops is invariant under $S$. Assume now that $k \geq 2$, fix an arbitrary cycle $C$ of length $k$ and define $S$ to map any non-loop arc to~$C$. For each loop, choose its image arbitrarily and extend $S$ to the whole set $\widehat{M}_n$ by additivity. Then the image of any graph without loops is $C$, so $S$ preserves the value $k$ of the cycle length. Thus, if we aim to classify additive mappings~$S$ that preserve a single value $k \geq 2$ of the cycle length, then it is natural to restrict our attention to the case where the set of non-loop arcs is invariant under~$S$. Under this condition, a mapping $S$ preserves the value~$2$ of the cycle length if and only if $S$ maps mutually inverse arcs to mutually inverse arcs. Therefore, the most interesting case is $k \geq 3$ which is precisely the case that is discussed in this paper.

Our paper is organized as follows. In Section 2 we consider additive mappings of unweighted directed graphs on $n$ vertices that preserve the values $2$ and $3$ of the cycle length. Lemma~\ref{lem_vertices_permutation} shows that if $n \geq 3$, then the set of non-loop arcs is invariant under such a mapping, and for any $n \geq 2$ the invariance of this set implies that the mapping acts bijectively on it and has the so-called canonical form. Namely, such a mapping corresponds to a renumbering of the vertices of the graph and, possibly, reversing the orientation of all its arcs. In Sections~3 and~4 we study additive mappings that preserve a single value $k \geq 3$ of the cycle length, with the set of non-loop arcs being again invariant. The main result of the paper is Theorem~\ref{th_vertices_permutation} which states that for $n \geq \max(k+1,5)$ the restriction of such a mapping to the set of non-loop arcs is bijective, and for any $n\geq \max(k, 4)$ the bijectivity of this restriction implies the canonical form of the mapping. Examples \ref{ex_prsrv_3_not_prsrv_2}--\ref{ex_prsrv_3_not_bijective2}, \ref{ex_prsrv_k_not_bijective}, and~\ref{ex_prsrv_2_not_prsrv_3} show that all the constraints of this theorem are essential.

\section{Additive mappings of graphs preserving the values $2$ and~$3$ of the cycle length}

For convenience, we introduce the following notation. Denote by $\mathcal{\widehat E}_n$ the set of graphs from $\widehat M_n$ with a single arc, by $\mathcal{\widehat D}_n$ its subset consisting of the graphs with a single loop, and by $\mathcal{\widehat N}_n$ the subset of graphs with a single non-loop arc. We identify the graph with a single arc from vertex $i$ to vertex $j$ with this arc and denote it by $e_{ij}$.

For a non-loop arc $e = e_{ij}$ we set $\overline{e} = e_{ji}$; then the arcs $e$ and $\overline{e}$ are called \textit{mutually inverse}. We call arcs $e_1, e_2 \in \mathcal{\widehat N}_n$ \textit{consecutive} if the terminal vertex of one of them coincides with the initial vertex of the other, and $e_1, e_2$ are not mutually inverse. We call distinct arcs $e_1, e_2 \in \mathcal{\widehat N}_n$ \textit{adjacent} if they are incident to at least one common vertex.

The proofs of Proposition~\ref{prop_adj_edges} and Lemma~\ref{lem_vertices_permutation} are based on the proofs of Proposition~2.8 and Lemma~2.10 in~\cite{GKV2}, but extend them to a more general setting.

\begin{proposition}\label{prop_adj_edges}
    Let $S:\widehat M_n\to \widehat M_n$ be an additive mapping of unweighted graphs that preserves the values $2$ and $3$ of the cycle length, and let the set $\mathcal{\widehat N}_n$ be invariant under the mapping~$S$.
    Then $S$ maps any two adjacent non-loop arcs to adjacent non-loop arcs; moreover, 
    \begin{enumerate}
        \item mutually inverse arcs are mapped to mutually inverse arcs,
        \item consecutive arcs are mapped to consecutive arcs,
        \item adjacent arcs that are neither consecutive nor mutually inverse are mapped to adjacent arcs that are neither consecutive nor mutually inverse.
    \end{enumerate}
\end{proposition}

\begin{proof}
\begin{enumerate}
    \item Since $S$ maps a cycle of length~$2$ to a cycle of length~$2$, for any arc $e \in \mathcal{\widehat N}_n$ we have $S(\overline{e}) = \overline{S(e)}$.
    \item Now let the arcs $e_1, e_2 \in\mathcal{\widehat{N}}_n$ be consecutive. Then there exists a cycle of length~$3$ of the form $e_1 \cup e_2 \cup e_3$. Since $S$ maps a cycle of length~$3$ to a cycle of length~$3$, $S(e_1 \cup e_2 \cup e_3)=S(e_1)\cup S(e_2)\cup S(e_3)$ is also a cycle of length~$3$. Hence $S(e_1)$ and $S(e_2)$ are consecutive arcs.
    \item Finally, let $e_1, e_2 \in \mathcal{\widehat{N}}_n$ be adjacent non-consecutive arcs that do not form a cycle of length $2$. Then the arcs $e_1$ and $\overline{e_2}$ are consecutive, so, by what has been proved, their images $S(e_1)$ and $S(\overline{e_2}) = \overline{S(e_2)}$ are also consecutive. Therefore $S(e_1)$ and $S(e_2)$ are adjacent non-consecutive arcs that do not form a cycle of length $2$. \qedhere
\end{enumerate}
\end{proof}

\begin{proposition}\label{prop_preserve_4_cycle_length}
    Let $S:\widehat M_n\to \widehat M_n$ be an additive mapping of unweighted graphs that preserves the values $2$ and $3$ of the cycle length, and let the set $\mathcal{\widehat N}_n$ be invariant under the mapping~$S$.
    Then $S$ preserves the value~$4$ of the cycle length.
\end{proposition}

\begin{proof}
    Consider an arbitrary cycle of length $4$ consisting of the arcs $e_{ij}, e_{jk}$, $e_{kl}, e_{li}$, where $i,j,k,l$ are pairwise distinct. Since $S$ preserves the value~$3$ of the cycle length, the arcs $S(e_{ij}), S(e_{jk}), S(e_{ki})$ form a cycle of length $3$. Two cases are possible: 
    \begin{enumerate}
        \item $S(e_{ij}) = e_{i'j'}$, $S(e_{jk}) = e_{j'k'}$, $S(e_{ki}) = e_{k'i'}$;
        \item $S(e_{ij}) = e_{j'i'}$, $S(e_{jk}) = e_{k'j'}$, $S(e_{ki}) = e_{i'k'}$,
    \end{enumerate}
    where $i',j',k'$ are pairwise distinct.
    We prove the statement for case~(1); case (2) is treated similarly. 
    
    Since $S$ maps cycles of length $2$ to cycles of length $2$, we have $S(e_{ik}) = e_{i'k'}$. Moreover, the arcs $S(e_{ik}), S(e_{kl}), S(e_{li})$ form a cycle of length $3$. Then we can write either $S(e_{kl}) = e_{k'l'}$, $S(e_{li}) = e_{l'i'}$, or $S(e_{kl}) = e_{l'i'}$, $S(e_{li}) = e_{k'l'}$, where $i', k', l'$ are pairwise distinct. The second option is impossible, since in this case the arcs $S(e_{ij}) = e_{i'j'}$ and $S(e_{li}) = e_{k'l'}$ are not consecutive, which contradicts Proposition~\ref{prop_adj_edges}. Thus the arcs $e_{ij}, e_{jk}$, $e_{kl}, e_{li}$ are mapped to $e_{i'j'}, e_{j'k'}, e_{k'l'}, e_{l'i'}$, respectively. The arcs $e_{jk}$ and $e_{kl}$ are consecutive, so $e_{j'k'}$ and $e_{k'l'}$ are also consecutive; in particular, $j' \neq l'$. Hence the mapping $S$ preserves the value~$4$ of the cycle length.
\end{proof}

\begin{lemma}\label{lem_vertices_permutation}
Let $S:\widehat M_n\to \widehat M_n$ be an additive mapping of unweighted graphs that preserves the values $2$ and $3$ of the cycle length, $n \geq 2$.
\begin{enumerate}
    \item If $n \geq 3$, then the set $\mathcal{\widehat N}_n$ is invariant under the mapping~$S$.
    \item If the set $\mathcal{\widehat N}_n$ is invariant under the mapping~$S$, then the restriction $S|_{\mathcal{\widehat N}_n}$ is bijective, and there exists a permutation $\tau\in S_n$ such that exactly one of the following two conditions holds simultaneously for all $i\ne j$:
    \begin{equation} \label{equation:canonical-form}
        S(e_{ij})=e_{\tau(i)\tau(j)} \quad \text{or} \quad S(e_{ij})=e_{\tau(j)\tau(i)}.
    \end{equation}
    In other words, if a graph $\widehat G$ contains no loops, then $S(\widehat G)$ is obtained from $\widehat G$ by a renumbering of the vertices and, possibly, reversing the orientation of all arcs.
\end{enumerate} 
\end{lemma}

\begin{proof}
Assume that $n \geq 3$ and consider an arbitrary arc $e \in \mathcal{\widehat{N}}_n$. Then there exist cycles $e \cup \overline{e}$ and $e \cup e_2 \cup e_3$ of lengths $2$ and $3$, respectively. Their images $S(e) \cup S(\overline{e})$ and $S(e) \cup S(e_2) \cup S(e_3)$ are also cycles of the same lengths. Reversing the orientation of all arcs in $e \cup e_2 \cup e_3$, we obtain a cycle $\overline{e} \cup \overline{e_2} \cup \overline{e_3}$ of length $3$ whose image $S(\overline{e}) \cup S(\overline{e_2}) \cup S(\overline{e_3})$ is again a cycle of length~$3$. Hence $S(e), S(\overline{e}) \in \mathcal{\widehat{N}}_n$, and item~(1) is proved.

We now prove~(2). Since the set $\mathcal{\widehat N}_n$ is invariant under the mapping~$S$, the restriction $S|_{\mathcal{\widehat N}_n}$ is well defined. The set $\mathcal{\widehat N}_n$ is finite, so, in order to prove that $S|_{\mathcal{\widehat N}_n}$ is bijective, it suffices to prove its injectivity.

If two arcs are adjacent, then, by Proposition~\ref{prop_adj_edges}, their images are also adjacent and therefore distinct. It remains to consider arcs $e_{ij}$ and~$e_{kl}$, where $i,j,k,l$ are pairwise distinct. By Proposition~\ref{prop_preserve_4_cycle_length}, $S$ preserves the value $4$ of the cycle length, so the image of the cycle $e_{ij} \cup e_{jk} \cup e_{kl} \cup e_{li}$ is also a cycle of length~$4$. Hence $S(e_{ij}) \neq S(e_{kl})$, and the bijectivity of $S|_{\mathcal{\widehat N}_n}$ is proved.

Next, we prove that $S$ has the canonical form~\eqref{equation:canonical-form}. By the additivity of the mapping~$S$, it suffices to consider only graphs that consist of a single non-loop arc. Thus, in view of Proposition~\ref{prop_adj_edges}, we need to show that one of the following two alternatives holds:

\begin{enumerate}[label=(\Roman*)]
    \item For any arcs $e_1, e_2\in\mathcal{\widehat N}_n$ the following statements hold:
    \begin{enumerate}[font=\upshape,label=(\alph*)]
        \item if the terminal vertex of $e_1$ coincides with the initial vertex of $e_2$, then the terminal vertex of $S(e_1)$ coincides with the initial vertex of $S(e_2)$;
        \item if $e_1$ and $e_2$ have the same initial vertex, then $S(e_1)$ and $S(e_2)$ also have the same initial vertex;
        \item if $e_1$ and $e_2$ have the same terminal vertex, then $S(e_1)$ and $S(e_2)$ also have the same terminal vertex.
    \end{enumerate}
    
    \item For any arcs $e_1, e_2\in\mathcal{\widehat N}_n$ the following statements hold:
    \begin{enumerate}[font=\upshape,label=(\alph*)]
        \item if the terminal vertex of $e_1$ coincides with the initial vertex of $e_2$, then the initial vertex of $S(e_1)$ coincides with the terminal vertex of $S(e_2)$;
        \item if $e_1$ and $e_2$ have the same initial vertex, then $S(e_1)$ and $S(e_2)$ have the same terminal vertex;
        \item if $e_1$ and $e_2$ have the same terminal vertex, then $S(e_1)$ and $S(e_2)$ have the same initial vertex.
    \end{enumerate}
\end{enumerate}

Indeed, in these cases we can define a permutation $\tau \in S_n$ as follows. If item (I) is satisfied, then, for any $1 \leq i \neq j \leq n$, $\tau(i)$ is the initial vertex of $S(e_{ij})$ which is well defined by~(Ib). Similarly, if item (II) is satisfied, then $\tau(i)$ is the terminal vertex of $S(e_{ij})$ which is well defined by~(IIb).

If $n=2$, then $e_{12}$ and~$e_{21}$ are the only two non-loop arcs, so (I) and~(II) follow from the fact that the mapping $S$ preserves cycles of length $2$. From now on we assume that $n\geq 3$.

Consider an arbitrary set of pairwise distinct indices $1\leq i,j,k\leq n$. Then the pair of consecutive arcs $e_{ij}, e_{jk}$ satisfies either (Ia) or (IIa), since the image of a cycle of length~$3$ containing these two arcs is again a cycle of length~$3$. Without loss of generality we may assume that (Ia) holds. First we show that for any index $1\leq t\leq n$ such that $j,k,t$ are pairwise distinct, the arcs $e_{jk}$ and $e_{kt}$ also satisfy condition (Ia). 
Two cases are possible:
\begin{enumerate}[label=(\roman*)]
    \item If $t=i$, then consider the cycle consisting of the arcs $e_{ij}, e_{jk}, e_{kt}$. Since the mapping~$S$ preserves cycles of length $3$ and the arcs $e_{ij}$ and $e_{jk}$ satisfy condition (Ia), any two arcs in this cycle satisfy (Ia).

    \item If $t\notin\{i,j,k\}$, then consider the cycle consisting of the arcs $e_{ij}, e_{jk},$ $e_{kt}, e_{ti}$. By Proposition~\ref{prop_preserve_4_cycle_length}, the mapping $S$ preserves cycles of length $4$, and by Proposition~\ref{prop_adj_edges}, $S$ maps consecutive arcs to consecutive arcs. Moreover, the arcs $e_{ij}$ and $e_{jk}$ satisfy condition (Ia); hence any two consecutive arcs in this cycle satisfy (Ia). In particular, this holds for the pair of arcs $e_{jk}$ and~$e_{kt}$.
\end{enumerate}

Similarly, for any index $1\leq s\leq n$ such that $i, j, s$ are pairwise distinct, the arcs $e_{si}$ and $e_{ij}$ also satisfy condition~(Ia). 

Now we prove that for any indices $1\leq p,q,r\leq n$, $q \notin \{ p, r\}$, the arcs $e_{pq}$ and $e_{qr}$ also satisfy condition (Ia). If $p = r$, then the statement follows from the fact that the mapping~$S$ preserves cycles of length~$2$. From now on we assume that $p, q, r$ are pairwise distinct. The following cases are possible:

\begin{enumerate}[label=(\roman*)]
    \item If $p = j$, $q = i$, then, by what has been proved above, $e_{ri}$ and $e_{ij}$ satisfy condition (Ia). Since $S$ preserves cycles of length~$2$, the arcs $e_{ji} = e_{pq}$ and $e_{ir} = e_{qr}$ also satisfy condition~(Ia).
    
    \item If $p = j$, $q \neq i$, then the following pairs of arcs satisfy condition (Ia):
    \[
    e_{ij}\ \text{and}\ e_{jk} \quad \Rightarrow \quad e_{ki}\ \text{and}\ e_{ip} \quad \Rightarrow \quad e_{ip}\ \text{and}\ e_{pq} \quad \Rightarrow \quad e_{pq}\ \text{and}\ e_{qr}.
    \]
    
    \item Let $p = k$, $q = j$. Since $e_{ij} = e_{iq}$ and $e_{jk} = e_{qp}$ satisfy condition (Ia) and $S$ preserves cycles of length~$2$, condition (Ia) is satisfied by
    \begin{align*}
    e_{pq}\ \text{and}\ e_{qi} \quad \Rightarrow \quad e_{ip}\ \text{and}\ e_{pq} \quad \Rightarrow \quad e_{pq}\ \text{and}\ e_{qr}.
    \end{align*}

    \item If $p = k$, $q \neq j$, then condition (Ia) is satisfied by
    \[
    e_{ij}\ \text{and}\ e_{jk} \quad \Rightarrow \quad e_{jk}\ \text{and}\ e_{pq} \quad \Rightarrow \quad e_{pq}\ \text{and}\ e_{qr}.
    \]

    \item If $p \neq j$ and $q = k$, then, by what has been proved above, $e_{jk} = e_{jq}$ and $e_{kp} = e_{qp}$ satisfy condition (Ia). Since $S$ preserves cycles of length~$2$, we obtain that condition (Ia) is satisfied by
    \[
    e_{pq}\ \text{and}\ e_{qj} \quad \Rightarrow \quad e_{jp}\ \text{and}\ e_{pq} \quad \Rightarrow \quad e_{pq}\ \text{and}\ e_{qr}.
    \]
    
    \item If $p\notin \{ j, k\}$ and $q \neq k$, then condition (Ia) is satisfied by
    \[
    e_{ij}\ \text{and}\ e_{jk} \quad \Rightarrow \quad e_{jk}\ \text{and}\ e_{kp} \quad \Rightarrow \quad e_{kp}\ \text{and}\ e_{pq} \quad \Rightarrow \quad e_{pq}\ \text{and}\ e_{qr}.
    \]
\end{enumerate}

As $S$ preserves cycles of length $2$, items (b) and (c) follow from item (a). Therefore item (I) holds for any pair of distinct arcs from~$\mathcal{\widehat N}_n$. The case where the arcs $e_{ij}$ and $e_{jk}$ satisfy condition (IIa) is completely similar. The proof is complete.
\end{proof}

\begin{corollary}
    Let $S:\widehat M_n\to \widehat M_n$ be an additive mapping of unweighted graphs that preserves the values $2$ and $3$ of the cycle length. Then
    $S$ preserves all values $k\geq 2$ of the cycle length.
\end{corollary}

\begin{proof}
    If $n = 2$, then the only possible value $k \geq 2$ of the cycle length is $k = 2$, and the statement follows from the hypothesis. If $n \geq 3$, then, by Lemma~\ref{lem_vertices_permutation}, the set $\mathcal{\widehat N}_n$ is invariant under the mapping~$S$, and, for an arbitrary graph $\widehat G$ without loops, its image $S(\widehat G)$ is obtained from $\widehat G$ by a renumbering of the vertices and, possibly, reversing the orientation of all arcs. Hence the image of any cycle of length $k \geq 2$ is again a cycle of length~$k$.
\end{proof}

\section{Additive mappings of graphs preserving the value $3$ of the cycle length}

\begin{proposition}\label{prop_bij_action_on_k_cycles}
    Let $S:\widehat M_n\to \widehat M_n$ be an additive mapping of unweighted graphs such that the restriction $S|_{\mathcal{\widehat N}_n}$ is bijective, and let $2\leq k\leq n$. The mapping $S$ preserves the value~$k$ of the cycle length if and only if it maps the set of cycles of length~$k$ bijectively onto itself.
\end{proposition}

\begin{proof}
    If $S$ preserves the value $k$ of the cycle length, then it maps the set of cycles of length~$k$ injectively into itself. This set is finite, therefore, such a mapping is bijective.

    The converse statement follows immediately from Definition~\ref{def_cycle_length_preserving}.
\end{proof}

\begin{lemma}\label{lemm_preserving_for_3}
    Let $n\geq4$, and let $S:\widehat M_n\to \widehat M_n$ be an additive mapping of unweighted graphs that preserves the value $3$ of the cycle length. If the set $\mathcal{\widehat N}_n$ is invariant under the mapping~$S$ and the restriction $S|_{\mathcal{\widehat N}_n}$ is bijective, then $S$ preserves the value~$2$ of the cycle length.
\end{lemma}

\begin{proof}
    Let $e_1, e_2\in\widehat{\mathcal{N}}_n$ be two arbitrary distinct non-loop arcs. Denote by $R_3(e_1, e_2)$ the number of arcs~$w \notin \{ e_1, e_2 \}$ such that $e_1, w$ belong to some cycle of length~$3$ and $e_2, w$ also belong to some cycle of length~$3$.

    Since $S$ preserves the value $3$ of the cycle length, if the arcs $e_1, w$ belong to a cycle of length~$3$, then $S(e_1), S(w)$ also belong to a cycle of length~$3$. Therefore,
    \[
        R_3(e_1, e_2) \leq R_3(S(e_1), S(e_2)).
    \]
    By our assumption, the restriction of $S$ to the set of non-loop arcs is bijective, so, by additivity, $S$ is bijective on the set of graphs without loops. Besides, by Proposition~\ref{prop_bij_action_on_k_cycles}, $S$ permutes the set of cycles of length $3$. Hence there exists an inverse mapping on the set of graphs without loops which also preserves cycles of length $3$, and we obtain the reverse inequality. Thus,
    \[
    R_3(e_1, e_2) = R_3(S(e_1), S(e_2)),
    \]
    that is, $S$ preserves $R_3(e_1, e_2)$ for any pair of arcs $e_1, e_2$.
    
    Next, we consider the possible values of $R_3(e_1, e_2)$ depending on the relative position of the arcs $e_1 \neq e_2$. Let $e_1 = e_{ij}$, $e_2 = e_{pq}$ for some $i\neq j$, $p\neq q$, $1\leq i,j,p,q\leq n$, and let~$w$ be an arc such that $e_1, w$ belong to the same cycle of length $3$ and $e_2, w$ belong to the same cycle of length $3$. The first condition on~$w$ implies that $w = e_{jt}$ or $w = e_{ti}$ for some index $t\notin \{i, j\}$, $1\leq t\leq n$. The following cases are possible:
    \begin{enumerate}
        \item $p=j$, $q=i$,  that is, $e_1 = \overline{e_2}$. Without loss of generality we may assume that $w = e_{jt}$. Since $e_2, w$ belong to the same cycle of length 3, it follows that $w = e_{rj}$ or $w = e_{ir}$ for some $r\notin \{i, j\}$, $1\leq r\leq n$. These alternatives imply that $t=j$ or $i=j$, respectively, so both of them are impossible. Therefore, no such arc~$w$ exists, and $R_3(e_1, e_2) = 0$.
        \item $p=j$, $q\neq i$. Then the only suitable arc $w$ is $w = e_{qi}$, and $R_3(e_1, e_2) = 1$.
        \item $p\neq j$, $q=i$. Similarly to case (II), we obtain $R_3(e_1, e_2) = 1$.
        \item $p=i$, $q\neq j$. Then the only suitable arcs $w$ are the arcs of the form $w = e_{ti}$, $t\notin\{i, j, q\}$, and $R_3(e_1, e_2) = n-3$.
        \item $p\neq i$, $q=j$. Then the only suitable arcs $w$ are the arcs of the form $w = e_{jt}$, $t\notin\{i, j, p\}$, and $R_3(e_1, e_2) = n-3$.
        \item The indices $i,j,p,q$ are pairwise distinct. Then only $w = e_{jp}$ and $w = e_{qi}$ are suitable; therefore, $R_3(e_1, e_2) = 2$.
    \end{enumerate}

    Note that for $n\geq 4$ the equality $e_1 = \overline{e_2}$ is equivalent to $R_3(e_1, e_2) = 0$. Since $S$ preserves the invariant $R_3$, this condition can be rewritten as $R_3(S(e_1), S(e_2)) = 0$, which is equivalent to $S(e_1) = \overline{S(e_2)}$. Therefore, the mapping $S$ preserves the value $2$ of the cycle length.
\end{proof}

The following examples show that the restriction $n \geq 4$ and the bijectivity of $S|_{\mathcal{\widehat N}_n}$ are essential in Lemma~\ref{lemm_preserving_for_3}.

\begin{example}\label{ex_prsrv_3_not_prsrv_2}
    Let $S:\widehat M_3\to \widehat M_3$ be defined as follows:
    \begin{align*}
    S(e_{12}) = e_{23},& \quad S(e_{23}) = e_{31}, \quad S(e_{31}) = e_{12},\\
    S(e_{ij}) = e_{ij},& \quad (i,j)\notin\{(1,2), (2,3), (3,1)\},
    \end{align*}
    and extended to the whole set $\widehat M_3$ by additivity, see Fig.~\ref{figure:example1}. Then the restriction $S|_{\mathcal{\widehat N}_3}$ is bijective and $S$ preserves the value $3$ of the cycle length, but $S(e_{12}\cup e_{21}) = e_{23}\cup e_{21}$, and $S$ does not preserve the value $2$ of the cycle length.

    \begin{figure}[ht]
    \centering
    \begin{tikzpicture}[
    vertex/.style={circle, draw=black},
    edge label/.style={midway,fill=white,inner sep=1pt,font=\small}
    ]

    \begin{scope}[xshift=0cm]
    \node[vertex] (v1) at (0,2) {1};
    \node[vertex] (v3) at (-1.732,0) {3};
    \node[vertex] (v2) at (1.732,0) {2};
    
    \draw[->,thick] (v1) to[loop above,min distance=10pt] node[edge label] {$e_{11}$} (v1);
    \draw[->,thick] (v2) to[loop right,min distance=10pt] node[edge label] {$e_{22}$} (v2);
    \draw[->,thick] (v3) to[loop left,min distance=10pt] node[edge label] {$e_{33}$} (v3);
    
    \draw[->,thick] (v1) to[bend left=15] node[edge label,right] {$e_{12}$} (v2);
    \draw[->,thick] (v2) to[bend left=15] node[edge label,left] {$e_{21}$} (v1);
    
    \draw[->,thick] (v2) to[bend left=15] node[edge label,below] {$e_{23}$} (v3);
    \draw[->,thick] (v3) to[bend left=15] node[edge label,below] {$e_{32}$} (v2);
    
    \draw[->,thick] (v3) to[bend left=15] node[edge label,left] {$e_{31}$} (v1);
    \draw[->,thick] (v1) to[bend left=15] node[edge label,right] {$e_{13}$} (v3);
\end{scope}

\draw[->,thick] (2.5,1) -- node[above] {$S$} (3.5,1);

\begin{scope}[xshift=7cm]
    \node[vertex] (v1) at (0,2) {1};
    \node[vertex] (v3) at (-1.732,0) {3};
    \node[vertex] (v2) at (1.732,0) {2};
    
    \draw[->,thick] (v1) to[loop above,min distance=10pt] node[edge label] {$S(e_{11})$} (v1);
    \draw[->,thick] (v2) to[loop right,min distance=10pt] node[edge label] {$S(e_{22})$} (v2);
    \draw[->,thick] (v3) to[loop left,min distance=10pt] node[edge label] {$S(e_{33})$} (v3);
    
    \draw[->,thick] (v1) to[bend left=15] node[edge label,right] {$S(e_{31})$} (v2);
    \draw[->,thick] (v2) to[bend left=15] node[edge label] {$S(e_{21})$} (v1);
    
    \draw[->,thick] (v2) to[bend left=15] node[edge label,below] {$S(e_{12})$} (v3);
    \draw[->,thick] (v3) to[bend left=15] node[edge label,below] {$S(e_{32})$} (v2);
    
    \draw[->,thick] (v3) to[bend left=15] node[edge label,left] {$S(e_{23})$} (v1);
    \draw[->,thick] (v1) to[bend left=15] node[edge label] {$S(e_{13})$} (v3);
\end{scope}

\end{tikzpicture}
\caption{A bijective mapping that preserves the value $3$ of the cycle length but does not preserve the value $2$. \label{figure:example1}}
\end{figure}
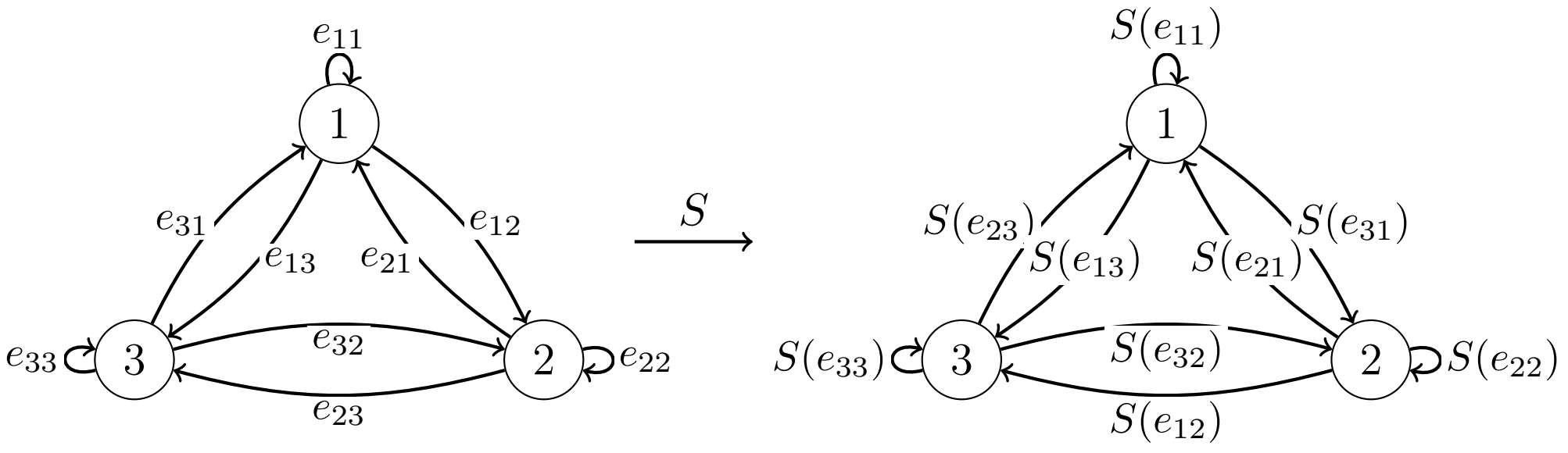
\end{example}

\begin{example} \label{ex_prsrv_3_not_bijective}
    Let $S:\widehat M_3\to \widehat M_3$ be defined as follows:
    \begin{align*}
    S(e_{21}) = e_{12},& \quad S(e_{13}) = e_{31}, \quad S(e_{32}) = e_{23},\\
    S(e_{ij}) = e_{ij},& \quad (i,j)\notin\{(2,1), (1,3), (3,2)\},
    \end{align*}
    and extended to the whole set $\widehat M_3$ by additivity, see Fig.~\ref{figure:example0}. Then $S$ preserves the value~$3$ of the cycle length and the set $\mathcal{\widehat N}_3$ is invariant under $S$, but the restriction $S|_{\mathcal{\widehat N}_3}$ is not bijective, and $S$ does not preserve the value~$2$ of the cycle length.

    \begin{figure}[t]
    \centering
    \begin{tikzpicture}[
    vertex/.style={circle, draw=black},
    edge label/.style={midway,fill=white,inner sep=1pt,font=\small}
    ]

    \begin{scope}[xshift=0cm]
    \node[vertex] (v1) at (0,2) {1};
    \node[vertex] (v3) at (-1.732,0) {3};
    \node[vertex] (v2) at (1.732,0) {2};
    
    \draw[->,thick] (v1) to[loop above,min distance=10pt] node[edge label] {$e_{11}$} (v1);
    \draw[->,thick] (v2) to[loop right,min distance=10pt] node[edge label] {$e_{22}$} (v2);
    \draw[->,thick] (v3) to[loop left,min distance=10pt] node[edge label] {$e_{33}$} (v3);
    
    \draw[->,thick] (v1) to[bend left=15] node[edge label,right] {$e_{12}$} (v2);
    \draw[->,thick] (v2) to[bend left=15] node[edge label,left] {$e_{21}$} (v1);
    
    \draw[->,thick] (v2) to[bend left=15] node[edge label,below] {$e_{23}$} (v3);
    \draw[->,thick] (v3) to[bend left=15] node[edge label,below] {$e_{32}$} (v2);
    
    \draw[->,thick] (v3) to[bend left=15] node[edge label,left] {$e_{31}$} (v1);
    \draw[->,thick] (v1) to[bend left=15] node[edge label,right] {$e_{13}$} (v3);
\end{scope}

\draw[->,thick] (2.5,1) -- node[above] {$S$} (3.5,1);

\begin{scope}[xshift=7cm]
    \node[vertex] (v1) at (0,2) {1};
    \node[vertex] (v3) at (-1.732,0) {3};
    \node[vertex] (v2) at (1.732,0) {2};
    
    \draw[->,thick] (v1) to[loop above,min distance=10pt] node[edge label] {$S(e_{11})$} (v1);
    \draw[->,thick] (v2) to[loop right,min distance=10pt] node[edge label] {$S(e_{22})$} (v2);
    \draw[->,thick] (v3) to[loop left,min distance=10pt] node[edge label] {$S(e_{33})$} (v3);
    
    \draw[->,thick] (v1) to[bend left=15] node[edge label,near end] {$S(e_{12}) = S(e_{21})$} (v2);
    
    \draw[->,thick] (v2) to[bend left=15] node[edge label,below] {$S(e_{23}) = S(e_{32})$} (v3);
    
    \draw[->,thick] (v3) to[bend left=15] node[edge label,near start] {$S(e_{31}) = S(e_{13})$} (v1);
\end{scope}

\end{tikzpicture}
\caption{A non-bijective mapping that preserves the value $3$ of the cycle length. \label{figure:example0}}
\end{figure}
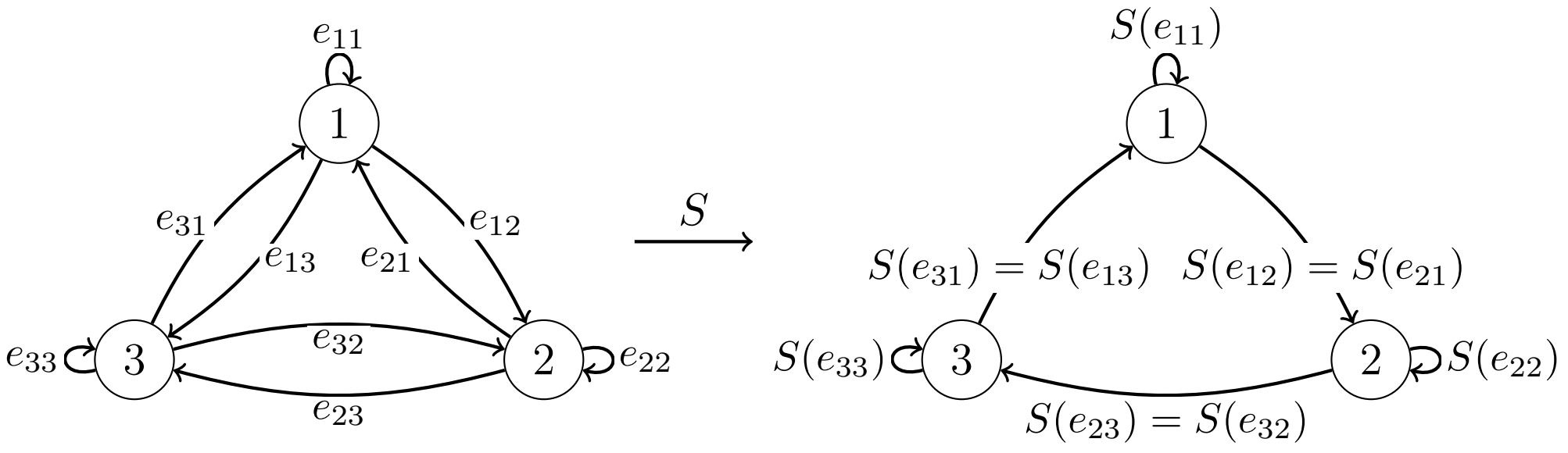
\end{example}

\begin{example} \label{ex_prsrv_3_not_bijective2}
    Let $S:\widehat M_4\to \widehat M_4$ be defined as follows:
    \begin{gather*}
    S(e_{12}) = S(e_{21}) = S(e_{34}) = S(e_{43}) = e_{12},\\
    S(e_{23}) = S(e_{32}) = S(e_{14}) = S(e_{41}) = e_{23},\\
    S(e_{13}) = S(e_{31}) = S(e_{24}) = S(e_{42}) = e_{31},\\
    S(e_{ii}) = e_{ii}, \quad 1 \leq i \leq 4,
    \end{gather*}
    and extended to the whole set $\widehat M_4$ by additivity. Then $S$ preserves the value $3$ of the cycle length and the set $\mathcal{\widehat N}_4$ is invariant under~$S$, but the restriction $S|_{\mathcal{\widehat N}_4}$ is not bijective, and $S$ does not preserve the value~$2$ of the cycle length.
\end{example}

\begin{lemma}\label{lemm_preserving_for_3_bijective}
    Let $S:\widehat M_n\to \widehat M_n$ be an additive mapping of unweighted graphs that preserves the value $3$ of the cycle length, and let the set $\mathcal{\widehat N}_n$ be invariant under the mapping~$S$. If $n \geq 5$, then the restriction $S|_{\mathcal{\widehat N}_n}$ is bijective.
\end{lemma}

\begin{proof}
    Since the set $\mathcal{\widehat N}_n$ is finite, it suffices to prove the injectivity of $S|_{\mathcal{\widehat N}_n}$. Suppose that $S(e_1) = S(e_2)$ for some arcs $e_1 \neq e_2$. If the arcs $e_1$ and $e_2$ are consecutive, then they are contained in some cycle of length $3$. Its image contains at most two distinct arcs and hence is not a cycle of length~$3$. This contradicts the fact that $S$ preserves the value $3$ of the cycle length.
    
    Therefore, the arcs $e_1 = e_{ij}$ and $e_2 = e_{pq}$ are not consecutive, $i \neq j$, $p \neq q$. Let $S(e_{ij}) = S(e_{pq}) = e_{ab}$. The following cases are possible:
    \begin{enumerate}
        \item The indices $i, j, p, q$ are pairwise distinct. Since $n \geq 5$, there exists $t \notin \{ i, j, p, q \}$. Considering the cycle $e_{ij} \cup e_{jt} \cup e_{ti}$, we obtain $S(e_{jt}) \in \{ e_{bc}, e_{ca} \}$, where $a, b, c$ are pairwise distinct. Similarly, the cycle $e_{pq} \cup e_{qt} \cup e_{tp}$ implies $S(e_{tp}) \in \{ e_{bd}, e_{da} \}$, where $a, b, d$ are pairwise distinct. Further, the image of the cycle $e_{tp} \cup e_{pj} \cup e_{jt}$ is also a cycle of length~$3$; hence $c = d$ and $S(e_{pj}) = e_{ab}$. Without loss of generality we may assume that $S(e_{jt}) = S(e_{qt}) = e_{bc}$ and $S(e_{ti}) = S(e_{tp}) = e_{ca}$.
        
        Since the arcs $e_{ij}$ and $e_{jq}$ are consecutive, they belong to the same cycle of length~$3$; hence $S(e_{ij}) = e_{ab}$ and $S(e_{jq})$ are also consecutive. Similarly, the arcs $S(e_{jq})$ and $S(e_{qt}) = e_{bc}$ are consecutive. Therefore, $S(e_{jq}) = e_{ca}$. Considering similarly the pairs of consecutive arcs $S(e_{pi})$ and $S(e_{ij}) = e_{ab}$, and also $S(e_{tp}) = e_{ca}$ and $S(e_{pi})$, we obtain $S(e_{pi}) = e_{bc}$.
        
        The image of the cycle $e_{ij} \cup e_{jp} \cup e_{pi}$ is a cycle of length~$3$ and equals $e_{ab} \cup S(e_{jp}) \cup e_{bc}$, whence $S(e_{jp}) = e_{ca}$. The image of the cycle $e_{jp} \cup e_{pq} \cup e_{qj}$ is also a cycle of length~$3$ and equals $e_{ca} \cup e_{ab} \cup S(e_{qj})$, whence $S(e_{qj}) = e_{bc}$. Thus, for the pair of consecutive arcs $e_{qj}$ and $e_{jt}$ we have $S(e_{qj}) = S(e_{jt}) = e_{bc}$, a contradiction.

        \item $p=j$, $q=i$. There exist $ 1 \leq r, t \leq n$ such that $i,j,r,t$ are pairwise distinct. Considering the cycle $e_{ij} \cup e_{jt} \cup e_{ti}$, we obtain $S(e_{jt}) \in \{ e_{bc}, e_{ca} \}$, where $a, b, c$ are pairwise distinct. Further, the cycle $e_{rj} \cup e_{ji} \cup e_{ir}$ implies $S(e_{rj}) \in \{ e_{bd}, e_{da} \}$, where $a, b, d$ are pairwise distinct. Since the image of the cycle $e_{rj} \cup e_{jt} \cup e_{tr}$ is a cycle of length~$3$, we obtain $c = d$ and $S(e_{tr}) = e_{ab}$. Hence $S(e_{ij}) = S(e_{tr})$, which gives a contradiction by case~(1).
        \item $p = i$, $q \neq j$. Choose $t \notin \{ i, j, q \}$. Considering the cycle $e_{ij} \cup e_{jt} \cup e_{ti}$, we obtain $S(e_{jt}) \in \{ e_{bc}, e_{ca} \}$, where $a, b, c$ are pairwise distinct. Further, the cycle $e_{iq} \cup e_{qj} \cup e_{ji}$ implies $S(e_{qj}) \in \{ e_{bd}, e_{da} \}$, where $a, b, d$ are pairwise distinct. Since the image of the cycle $e_{qj} \cup e_{jt} \cup e_{tq}$ is a cycle of length~$3$, we obtain $c = d$ and $S(e_{tq}) = e_{ab}$. Hence $S(e_{ij}) = S(e_{tq})$, which gives a contradiction by case~(1).
        \item $p \neq i$, $q = j$. The proof is similar to the previous case. \qedhere
    \end{enumerate}
\end{proof}

\section{Additive mappings of graphs preserving the value $k \geq 4$ of the cycle length}

\begin{proposition}\label{prop_existence_k_cycle}
    Let $4\leq k\leq n$. Distinct arcs $u, v \in \mathcal{\widehat N}_n$ belong to the same cycle of length~$k$ if and only if the following three conditions hold:
    \begin{enumerate}[label=(\alph*)]
        \item the initial vertex of $u$ does not coincide with the initial vertex of $v$;
        \item the terminal vertex of $u$ does not coincide with the terminal vertex of $v$;
        \item the arcs $u$ and $v$ are not mutually inverse.
    \end{enumerate}
\end{proposition}

\begin{proof}
    Let $u = e_{ij}$ and $v = e_{pq}$. Conditions (a)--(c) are satisfied if and only if one of the following two alternatives holds: either $u$ and $v$ are not adjacent, that is, the indices $i, j, p, q$ are pairwise distinct, or $u$ and $v$ are consecutive arcs, that is, $j = p$ and $i \neq q$, or $i = q$ and $j \neq p$. One can easily check that a cycle of length~$k$ containing $u$ and~$v$ exists precisely in these two cases.
\end{proof}

\begin{lemma}\label{lemm_preserving_for_k}
    Let $S:\widehat M_n\to \widehat M_n$ be an additive mapping of unweighted graphs that preserves the value $k$ of the cycle length, $4\leq k\leq n$. If the set $\mathcal{\widehat N}_n$ is invariant under the mapping~$S$ and the restriction $S|_{\mathcal{\widehat N}_n}$ is bijective, then $S$ preserves the value~$3$ of the cycle length.
\end{lemma}

\begin{proof}
    Let $e_1, e_2\in\widehat {\mathcal{N}}_n$ be two arbitrary distinct non-loop arcs. Denote by $R_k(e_1, e_2)$ the number of arcs $w \notin \{ e_1, e_2 \}$ such that $e_1, w$ belong to some cycle of length $k$ and $e_2, w$ also belong to some cycle of length~$k$. Since $S$ preserves the value $k$ of the cycle length, it preserves $R_k(e_1, e_2)$ for any pair of arcs $e_1, e_2$. The proof of this fact is similar to that for the case $k = 3$ in Lemma~\ref{lemm_preserving_for_3}.

    We consider all possible values of $R_k(e_1, e_2)$ depending on the relative position of the arcs $e_1 \neq e_2$. Let $e_1 = e_{ij}$, $e_2 = e_{pq}$ for some $i\neq j$, $p\neq q$, $1\leq i,j,p,q\leq n$. The following cases are possible:
    \begin{enumerate}
        \item $p=j$, $q=i$. The suitable arcs $w$ are all arcs $e_{xy}$ such that $x\neq y$ and $x,y\notin \{i,j\}$. That is, $R_k(e_1, e_2) = (n-2)\cdot(n-3) = n^2-5n+6$.
        \item $p=j$, $q\neq i$. Then the suitable arcs $w$ are all arcs $e_{xy}$ such that $x\neq y$, $x\notin\{i, j\}$, $y\notin\{j, q\}$. We obtain $R_k(e_1, e_2) = (n-2)^2 - (n-3) = n^2 - 5n + 7$.
        \item $p\neq j$, $q=i$. Similarly to case (2), we have $R_k(e_1, e_2) = n^2 - 5n + 7$.
        \item $p=i$, $q\neq j$. In this case the suitable arcs $w$ are all arcs $e_{xy}$ such that $x\neq y$, $x\neq i$, $y\notin\{j, q\}$, except $e_{ji}$ and $e_{qi}$. Then $R_k(e_1, e_2) = (n-1)\cdot(n-2) - (n-3) - 2 = n^2-4n+3$.
        \item $p\neq i$, $q=j$. Similarly to case (4), we obtain $R_k(e_1, e_2) = n^2-4n+3$.
        \item The indices $i,j,p,q$ are pairwise distinct. Then the suitable arcs $w$ are all arcs $e_{xy}$ such that $x\neq y$, $x\notin\{i,p\}$, $y\notin\{j,q\}$, except $e_{ji}$ and $e_{qp}$. Therefore, we have $R_k(e_1, e_2) = (n-2)^2 - (n-4) - 2 = n^2 - 5n + 6$.
    \end{enumerate}

    To prove that the mapping $S$ preserves the value~$3$ of the cycle length, it suffices to show that if the arcs $e_1$ and $e_2$ are consecutive, then $S(e_1)$ and $S(e_2)$ are also consecutive.
    
    Let $e_1$ and $e_2$ be consecutive arcs. Then, by the above arguments, we have \[R_k(S(e_1), S(e_2)) = R_k(e_1, e_2) = n^2-5n+7.\]
    For $n\geq 5$ we immediately obtain that $S(e_1)$ and $S(e_2)$ are consecutive arcs. For $n=4$ there is one more possibility. Namely, $4^2-5\cdot4+7 = 3 = 4^2 - 4\cdot 4 + 3$, that is, $S(e_1)$ and $S(e_2)$ can be adjacent non-consecutive arcs. But the arcs $e_1$ and $e_2$ are contained in some cycle of length $k$, and the mapping $S$ preserves the value~$k$ of the cycle length. Therefore, $S(e_1)$ and $S(e_2)$ are also contained in some cycle of length $k$. This contradicts the fact that these arcs are adjacent non-consecutive.
\end{proof}

Let us show that the condition of bijectivity of $S|_{\mathcal{\widehat N}_n}$ is essential in Lemma~\ref{lemm_preserving_for_k}.

\begin{example} \label{ex_prsrv_k_not_bijective}
    Let $n = k \geq 4$, and let $S:\widehat M_n\to \widehat M_n$ be defined as follows:
    \begin{gather*}
    \begin{aligned}
    S(e_{12}) = S(e_{13}) &= \ldots = S(e_{1n}) = e_{12},\\
    S(e_{21}) = S(e_{23}) &= \ldots = S(e_{2n}) = e_{23},\\
    &\hphantom{\, = \:} \ldots\\
    S(e_{n1}) = S(e_{n2}) &= \ldots = S(e_{n,n-1}) = e_{n1},
    \end{aligned}\\
    S(e_{ii}) = e_{ii}, \quad 1 \leq i \leq n,
    \end{gather*}
    and extended to the whole set $\widehat M_n$ by additivity. Then $S$ preserves the value $n$ of the cycle length and the set $\mathcal{\widehat N}_n$ is invariant under~$S$, but the restriction $S|_{\mathcal{\widehat N}_n}$ is not bijective, and $S$ does not preserve the value~$3$ of the cycle length.
\end{example}

\begin{lemma}\label{lemm_preserving_for_k_bijective}
    Let $S:\widehat M_n\to \widehat M_n$ be an additive mapping of unweighted graphs that preserves the value $k$ of the cycle length, $k \geq 4$, and let the set $\mathcal{\widehat N}_n$ be invariant under the mapping~$S$. If $n \geq k + 1$, then the restriction $S|_{\mathcal{\widehat N}_n}$ is bijective.
\end{lemma}

\begin{proof}
    Note that if distinct arcs $e_1, e_2 \in \mathcal{\widehat N}_n$ belong to the same cycle of length $k$, then its image is also a cycle of length~$k$ and hence contains $k$ distinct arcs. Therefore, $S(e_1) \neq S(e_2)$.
    
    Moreover, if $S$ maps every pair of consecutive arcs to consecutive arcs, then $S$ preserves the value $3$ of the cycle length; hence $S|_{\mathcal{\widehat N}_n}$ is bijective by Lemma~\ref{lemm_preserving_for_3_bijective}.

    Now suppose that there exists a pair of consecutive non-loop arcs $e_{ij}$ and $e_{jl}$ such that the arcs $S(e_{ij})$ and $S(e_{jl})$ are not consecutive. Take a cycle $C$ of length $k$ containing $e_{ij}$ and~$e_{jl}$. Since $n \geq k + 1$, there exists a vertex $t$ that is not incident to any arc of the cycle $C$. Let $C'$ be the cycle obtained from $C$ by replacing the arcs $e_{ij}$ and $e_{jl}$ with $e_{it}$ and $e_{tl}$. Then $S(C)$ and $S(C')$ are cycles of length~$k$, and $S(C')$ is obtained from $S(C)$ by replacing the arcs $S(e_{ij})$ and $S(e_{jl})$ with the arcs $S(e_{it})$ and $S(e_{tl})$. This is possible only if $\{ S(e_{ij}), S(e_{jl}) \} = \{ S(e_{it}), S(e_{tl}) \}$. Since the arcs $e_{ij}$ and $e_{tl}$ are contained in the same cycle of length $k$, the option $S(e_{ij}) = S(e_{tl})$ is impossible. Therefore,
    \begin{equation} \label{equation:bijective-k-cycles}
    S(e_{ij}) = S(e_{it}) \quad \text{and} \quad S(e_{jl}) = S(e_{tl}).
    \end{equation}
    
    The cycle $C$ and the vertex $t$ were chosen arbitrarily, so Eq.~\eqref{equation:bijective-k-cycles} holds for all $t \notin \{ i, l \}$. Note that then for any $p \neq l$ the arcs $S(e_{jl})$ and $S(e_{lp})$ are consecutive. Indeed, otherwise choose $t \notin \{ i, j, l, p\}$ and $s \notin \{ j, l, p, t \}$, which is possible since $n \geq k + 1 \geq 5$. Then $S(e_{jl}) = S(e_{tl})$, and, applying the statement proved above to the pair of arcs $e_{tl}$ and $e_{lp}$, we obtain $S(e_{tl}) = S(e_{ts})$. Hence $S(e_{jl}) = S(e_{ts})$, but the arcs $e_{jl}$ and $e_{ts}$ are contained in the same cycle of length $k$, a contradiction.

    Similarly, in this case for any $q \neq i$ the arcs $S(e_{qi})$ and $S(e_{ij})$ are consecutive.

    Let us return to the cycle~$C$. Since its image $S(C)$ is a cycle of length~$k$, it contains two arcs consecutive with $S(e_{ij})$. Therefore, there exists an arc $e_{rs}$ of the cycle $C$ that is not consecutive with~$e_{ij}$ (that is, $i,j,r,s$ are pairwise distinct) such that the arcs $S(e_{ij})$ and $S(e_{rs})$ are consecutive. Without loss of generality we may assume that $S(e_{ij}) = e_{ab}$ and $S(e_{rs}) = e_{bc}$, where $a,b,c$ are pairwise distinct. Hence the image of at least one of the two arcs of the cycle $C$ consecutive with~$e_{rs}$ is not consecutive with $S(e_{rs})$. It follows from what has been proved above that then either $S(e_{sx})$ is consecutive with~$S(e_{rs})$ for any $x \neq s$, or $S(e_{yr})$ is consecutive with~$S(e_{rs})$ for any $y \neq r$. If $x \neq j$ and $y \neq i$, the vertex~$b$ cannot be a terminal vertex of such an arc, since $e_{ij}$ belongs to the same cycle of length $k$ with each of the arcs $e_{sx}$ and $e_{yr}$. Therefore, the initial vertex of this arc is~$c$.

    Consider the first case. Since the arc $e_{rs}$ is contained in the cycle~$C$ and is not consecutive with~$e_{ij}$, while the arc $e_{jl}$ is contained in the same cycle and is consecutive with~$e_{ij}$, we have $s \neq l$. Hence we can set $t = s$ and $x = l$; therefore, by what has been proved above, $S(e_{sl}) = S(e_{jl})$ leaves the vertex~$c$. The arc $S(e_{si})$ is consecutive with each of the arcs $S(e_{ij}) = e_{ab}$ and $S(e_{rs}) = e_{bc}$; hence $S(e_{si}) = e_{ca}$. We obtain that $S(e_{jl})$ and $S(e_{si})$ leave the same vertex, but, since $i,j,l,s$ are pairwise distinct, they must belong to the same cycle of length~$k$, a contradiction.

    Now we turn to the second case. Then the arc consecutive with~$e_{rs}$ whose image is not consecutive with $S(e_{rs})$ has the form $e_{sw}$; moreover, since $e_{rs} \neq e_{ij}$, we obtain $e_{sw} \neq e_{jl}$ and $w \neq l$. Therefore, $S(e_{rz}) = S(e_{rs})$ for all $z \notin \{ r, w \}$. If $r \neq l$, then we can set $t = r$ and $z = l$; hence the arcs $S(e_{ij})$ and $S(e_{jl}) = S(e_{rl}) = S(e_{rs})$ are consecutive, a contradiction.
    
    Thus we have $r = l$. For $y = j$ we obtain that $S(e_{jl})$ leaves the vertex $c$, that is, $S(e_{jl}) = e_{cd}$. Since the arcs $e_{ij}, e_{jl}, e_{ls}$ belong to the cycle~$C$ of length~$k \geq 4$, their images also belong to the same cycle of length~$k$; hence $a,b,c,d$ are pairwise distinct. Repeating the previous arguments, we find an arc $e_{uv}$ of the cycle~$C$ that is not consecutive with~$e_{jl}$ such that the arcs $S(e_{jl})$ and $S(e_{uv})$ are consecutive. Then $S(e_{uv})$ leaves the vertex~$d$. Similarly, we obtain that this is possible only if $v = i$ and $S(e_{ui}) = e_{da}$. Therefore, the image of the cycle $C$ contains the arcs $e_{ab}, e_{bc}, e_{cd}, e_{da}$, which form a cycle of length~$4$. Since the mapping~$S$ preserves the value~$k$ of the cycle length, we obtain $k = 4$ and $C = e_{ij} \cup e_{jl} \cup e_{ls} \cup e_{si}$.

    For all $t \notin \{ i, l \}$ we have $S(e_{it}) = e_{ab}$, $S(e_{lt}) = e_{bc}$, $S(e_{tl}) = e_{cd}$ and $S(e_{ti}) = e_{da}$. By what has been proved, the arc $S(e_{li})$ is consecutive with each of the arcs $S(e_{ij})$ and $S(e_{jl})$. Hence $S(e_{li}) \in \{ e_{bc}, e_{da} \}$. The image of the cycle $e_{li} \cup e_{ij} \cup e_{js} \cup e_{sl}$ equals $S(e_{li}) \cup e_{ab} \cup S(e_{js}) \cup e_{cd}$ and is a cycle of length~$4$; thus $S(e_{js}) \in \{ e_{da}, e_{bc} \}$. Choose $t \notin \{ i, j, l, s\}$; then the arc $e_{js}$ belongs to the same cycle of length~$4$ with each of the arcs $e_{lt}$ and~$e_{ti}$, and either $S(e_{js}) = S(e_{lt})$ or $S(e_{js}) = S(e_{ti})$, a contradiction.
\end{proof}

\begin{theorem}\label{th_vertices_permutation}
    Let $S:\widehat M_n\to \widehat M_n$ be an additive mapping of unweighted graphs that preserves the value $k\geq 3$ of the cycle length, $n\geq\max(k,4)$, and let the set $\mathcal{\widehat N}_n$ be invariant under the mapping~$S$.
    \begin{enumerate}
        \item If $n \geq \max(k+1,5)$, then the restriction $S|_{\mathcal{\widehat N}_n}$ is bijective.
        \item If the restriction $S|_{\mathcal{\widehat N}_n}$ is bijective, then there exists a permutation $\tau\in S_n$ such that exactly one of the following two conditions holds simultaneously for all $i\ne j$:
        \[
        S(e_{ij})=e_{\tau(i)\tau(j)} \quad \text{or} \quad S(e_{ij})=e_{\tau(j)\tau(i)}.
        \]
        In other words, if a graph $\widehat G$ contains no loops, then $S(\widehat G)$ is obtained from $\widehat G$ by a renumbering of the vertices and, possibly, reversing the orientation of all arcs.
        \end{enumerate}
\end{theorem}

\begin{proof}
    Item (1) follows immediately from Lemmas~\ref{lemm_preserving_for_3_bijective} and~\ref{lemm_preserving_for_k_bijective}.

    Let us prove item (2). If $k\geq 4$, then, by Lemma~\ref{lemm_preserving_for_k}, the mapping~$S$ preserves the value $3$ of the cycle length. For $k=3$ this holds by the hypothesis of the theorem. Thus, in both cases we can apply Lemma~\ref{lemm_preserving_for_3}. Therefore, the mapping $S$ preserves the value~$2$ of the cycle length. Further, applying Lemma~\ref{lem_vertices_permutation}, we obtain the required statement.
\end{proof}

The following example shows that preserving the value $2$ of the cycle length is not sufficient.

\begin{example}\label{ex_prsrv_2_not_prsrv_3}
    Let $n\geq 3$, and let the mapping $S: \widehat M_n\to\widehat M_n$ be defined as follows:
    \begin{gather*}
        S(e_{12}) = e_{21},\quad S(e_{21}) = e_{12},\\
        S(e_{ij}) = e_{ij},\quad (i,j)\notin\{(1,2), (2,1)\},
    \end{gather*}
    and extended to the whole set $\widehat M_n$ by additivity, see Fig.~\ref{figure:example2}. It is clear that $S$ preserves the value $2$ of the cycle length. At the same time, $S(e_{12}\cup e_{23}\cup e_{31}) = e_{21}\cup e_{23}\cup e_{31}$; therefore, the mapping $S$ does not preserve the value~$3$ of the cycle length and is not of the form described in Lemma~\ref{lem_vertices_permutation}.

        \begin{figure}[h]
        \centering   
        \begin{tikzpicture}[
        vertex/.style={circle, draw=black},
        edge label/.style={midway,fill=white,inner sep=1pt,font=\small}
        ]
    
        \begin{scope}[xshift=0cm]
        \node[vertex] (v1) at (0,2) {1};
        \node[vertex] (v3) at (-1.732,0) {3};
        \node[vertex] (v2) at (1.732,0) {2};
        
        \draw[->,thick] (v1) to[loop above,min distance=10pt] node[edge label] {$e_{11}$} (v1);
        \draw[->,thick] (v2) to[loop right,min distance=10pt] node[edge label] {$e_{22}$} (v2);
        \draw[->,thick] (v3) to[loop left,min distance=10pt] node[edge label] {$e_{33}$} (v3);
        
        \draw[->,thick] (v1) to[bend left=15] node[edge label,right] {$e_{12}$} (v2);
        \draw[->,thick] (v2) to[bend left=15] node[edge label,left] {$e_{21}$} (v1);
        
        \draw[->,thick] (v2) to[bend left=15] node[edge label,below] {$e_{23}$} (v3);
        \draw[->,thick] (v3) to[bend left=15] node[edge label,below] {$e_{32}$} (v2);
        
        \draw[->,thick] (v3) to[bend left=15] node[edge label,left] {$e_{31}$} (v1);
        \draw[->,thick] (v1) to[bend left=15] node[edge label,right] {$e_{13}$} (v3);
    \end{scope}
    
    \draw[->,thick] (2.5,1) -- node[above] {$S$} (3.5,1);
    
    \begin{scope}[xshift=7cm]
        \node[vertex] (v1) at (0,2) {1};
        \node[vertex] (v3) at (-1.732,0) {3};
        \node[vertex] (v2) at (1.732,0) {2};
        
        \draw[->,thick] (v1) to[loop above,min distance=10pt] node[edge label] {$S(e_{11})$} (v1);
        \draw[->,thick] (v2) to[loop right,min distance=10pt] node[edge label] {$S(e_{22})$} (v2);
        \draw[->,thick] (v3) to[loop left,min distance=10pt] node[edge label] {$S(e_{33})$} (v3);
        
        \draw[->,thick] (v1) to[bend left=15] node[edge label,right] {$S(e_{21})$} (v2);
        \draw[->,thick] (v2) to[bend left=15] node[edge label] {$S(e_{12})$} (v1);
        
        \draw[->,thick] (v2) to[bend left=15] node[edge label,below] {$S(e_{23})$} (v3);
        \draw[->,thick] (v3) to[bend left=15] node[edge label,below] {$S(e_{32})$} (v2);
        
        \draw[->,thick] (v3) to[bend left=15] node[edge label,left] {$S(e_{31})$} (v1);
        \draw[->,thick] (v1) to[bend left=15] node[edge label] {$S(e_{13})$} (v3);
    \end{scope}
    
    \end{tikzpicture}
    \caption{A bijective mapping that preserves the value $2$ of the cycle length but does not preserve the value $3$. \label{figure:example2}}
    \end{figure}
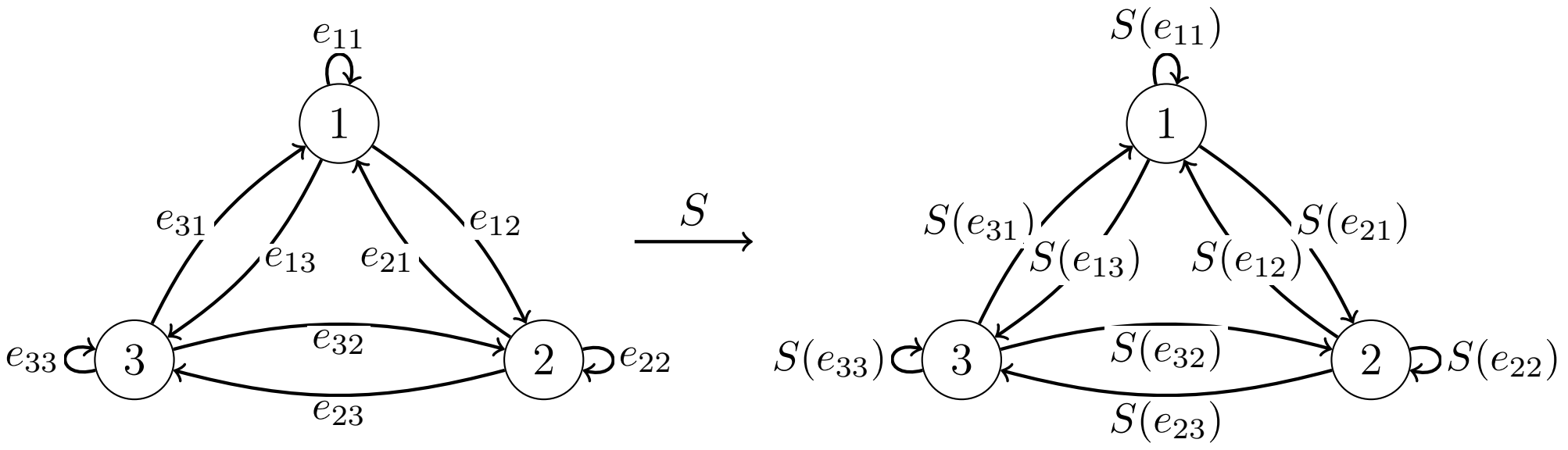
\end{example}

\bibliographystyle{abbrvurl}
\bibliography{thebibliography-arXiv}

@article{GKV2,
  author = {Vlasov, A. V. and
            Guterman, A. E. and
            Kreines, E. M.},
  title = {Linear transformations preserving the minimal values of the cyclicity index of tropical matrices},
  journal = {J. Math. Sci.},
  volume = {281},
  number = {2},
  year = {2024},
  pages = {234--245},
  doi = {10.1007/s10958-024-07096-y}
}

@article{BP1,
  author = {Beasley, L. B. and Pullman, N. J.},
  title = {Term-rank, permanent, and rook-polynomial preservers},
  journal = {Linear Algebra Appl.},
  volume = {90},
  year = {1987},
  pages = {33--46},
  doi = {10.1016/0024-3795(87)90302-8}
}

@article{BP2,
  author = {Beasley, L. B. and Pullman, N. J.},
  title = {Linear operators that strongly preserve primitivity},
  journal = {Linear Multilinear Algebra},
  volume = {25},
  year = {1989},
  pages = {205--213},
  doi = {10.1080/03081088908817942}
}

@article{Di,
  author = {Dieudonn{\'e}, J.},
  title = {Sur une g{\'e}n{\'e}ralisation du groupe orthogonal {\`a}
           quatre variables},
  journal = {Arch. Math.},
  volume = {1},
  year = {1949},
  pages = {282--287},
  doi = {10.1007/BF02038756}
}

@article{F,
  author = {Frobenius, G.},
  title = {{{\"U}ber die Darstellung der endlichen Gruppen durch
           lineare Substitutionen}},
  journal = {Sitzungsber. Deutsch. Akad. Wiss., Berlin},
  year = {1897},
  pages = {994--1015}
}

@article{GKT,
  author = {Guterman, A. and Kreines, E. and Thomassen, C.},
  title = {Linear transformations of tropical matrices preserving
           the cyclicity index},
  journal = {Special Matrices},
  volume = {9},
  year = {2021},
  pages = {112--118},
  doi = {10.1515/spma-2020-0128}
}

@article{GKV1,
  author = {Guterman, A. and Kreines, E. and Vlasov, A.},
  title = {Non-surjective linear transformations of tropical matrices
           preserving the cyclicity index},
  journal = {Kybernetika},
  volume = {58},
  number = {5},
  year = {2022},
  pages = {691--707},
  doi = {10.14736/kyb-2022-5-0691}
}

@article{GM1,
  author = {Guterman, A. and Maksaev, A.},
  title = {Maps preserving scrambling index},
  journal = {Linear Multilinear Algebra},
  volume = {66},
  number = {4},
  year = {2018},
  pages = {840--851},
  doi = {10.1080/03081087.2017.1329814}
}

@misc{GVZ2,
  title={Linear transformations of tropical matrices preserving individual values of the cyclicity index},
  author={Guterman, A. and Vlasov, A. and Zhilina, S.},
  year={2026},
  archivePrefix={arXiv},
  primaryClass={math.CO},
  eprint={2610.32897}
}

@book{HOW,
  author = {Heidergott, B. and Olsder, G. J. and van der Woude, J.},
  title = {Max Plus at Work},
  series = {Princeton Series in Applied Mathematics},
  publisher = {Princeton University Press},
  year = {2006},
  doi = {10.1515/9781400865239}
}

@article{Pierce,
  author = {Pierce, S. and others},
  title = {A survey of linear preserver problems},
  journal = {Linear Multilinear Algebra},
  volume = {33},
  year = {1992},
  pages = {1--129},
  doi = {10.1080/03081089208818176}
}

@article{Schur,
  author = {Schur, I.},
  title = {{Einige Bemerkungen zur Determinantentheorie}},
  journal = {Sitzungsber. Preuss. Akad. Wiss. Berlin},
  year = {1925},
  pages = {454--463}
}

@article{Thomassen,
  title     = {Whitney's $2$-switching theorem, cycle spaces, and arc mappings of directed graphs},
  author    = {Thomassen, Carsten},
  journal   = {J. Combin. Theory Ser. B},
  volume    = {46},
  number    = {3},
  pages     = {257--291},
  year      = {1989},
  doi       = {10.1016/0095-8956(89)90050-6}
}

@article{Woodall,
  title     = {Cycle lengths and circuit matroids of graphs},
  author    = {Woodall, D. R.},
  journal   = {Discrete Math.},
  volume    = {105},
  number    = {1--3},
  pages     = {269--273},
  year      = {1992},
  doi       = {10.1016/0012-365X(92)90149-A}
}

\end{document}